\documentclass[12pt]{article}
\usepackage[margin=2.5cm]{geometry}
\usepackage{graphicx,url}
\usepackage{amssymb,bm}
\usepackage[english]{babel}

\def\no{\noindent}
\newtheorem{defi}{Definition}

\newtheorem{rem}{Remark}
\newtheorem{theo}{Theorem}
\newtheorem{cor}{Corollary}
\newtheorem{lem}{Lemma}
\newtheorem{ass}{Assumption}
\def\proof{\no\underline{Proof.}~}
\def\QED{\mbox{$\square$}}
\def\no{\noindent}
\def\RR{\mathbb{R}}

\def\pmatrix{\left(\begin{array}}
\def\endpmatrix{\end{array}\right)}
\def\dd{{\mathrm d}}

\def\bfb{{\bf b}}
\def\bfc{{\bf c}}

\def\bfe{{\bf e}}

\def\aa{{\alpha}}

\def\eps{\varepsilon}

\def\bfgamma{\bm{\gamma}}
\def\balfa{\bm{\alpha}}
\def\bgamma{\bar{\gamma}}
\def\hgamma{\hat{\gamma}}

\def\II{{\cal I}}
\def\PP{{\cal P}}
\def\dd{{\mathrm{d}}}

\def\diag{{\rm diag}}

\def\heta{{\hat\eta}}
\def\hphi{{\hat\phi}}
\def\haa{{\hat\alpha}}
\def\hb{\hat{b}}
\def\hc{\hat{c}}
\def\hu{\hat{u}}

\title{Multiple invariants conserving Runge-Kutta type methods for Hamiltonian problems}

\author{Luigi Brugnano\footnote{Dipartimento di Matematica e Informatica ``U.\,Dini'', Universit\`a di Firenze, Italy
({\tt luigi.brugnano@unifi.it})}
\and Yajuan Sun\footnote{Academy of Mathematics and Systems Science, Chinese Academy of Sciences, Beijing, China ({\tt sunyj@lsec.cc.ac.cn})}}
\date{\em\small Warmly dedicated to celebrate the 80-th birthday of John Butcher}

\begin{document}
\maketitle
\begin{abstract}
In a recent series of papers, the class of energy-conserving Runge-Kutta methods named {\em Hamiltonian BVMs (HBVMs)} has been defined and studied. Such methods have been further generalized for the efficient solution of general conservative problems, thus providing the class of {\em Line Integral Methods (LIMs)}. In this paper we derive a further extension, which we name {\em Enhanced Line Integral Methods (ELIMs)}, more tailored for Hamiltonian problems, allowing for the conservation of multiple invariants of the continuous dynamical system. The analysis of the methods is fully carried out and some numerical tests are reported, in order to confirm the theoretical achievements.

\medskip
\no{\bf Keywords:}  Hamiltonian problems, energy-conserving methods, multiple invariants, discrete line-integral methods, HBVMs, LIMs, ELIMs, EHBVMs.

\medskip
\no{\bf MSC:} 65P10, 65L05.
\end{abstract}

\section{Introduction} Hamiltonian problems arise in many fields of application, ranging from the nano-scale of molecular dynamics to the macro-scale of celestial mechanics. Such problems are in the following form:
\begin{equation}\label{ham}
y' = J\nabla H(y), \qquad y(0) = y_0\in\RR^{2m},
\end{equation}
where the state vector is often partitioned as $$y = \pmatrix{c} q\\ p\endpmatrix, \qquad q,p\in\RR^m,$$ with $q$ the vector of the positions and $p$ the vector of the generalized momenta.
Moreover,
\begin{equation}\label{J}
J = \pmatrix{cc} 0 & I_m\\ -I_m &0\endpmatrix = -J^T = -J^{-1},
\end{equation}
and $H(y)\equiv H(q,p)$ is the {\em Hamiltonian} function defining the problem.  From (\ref{ham}) and (\ref{J}), it is straightforward to derive that $H(y(t))\equiv H(y_0)$ for $t\ge0$, since
$$\frac{\dd}{\dd t} H(y(t)) = \nabla H(y(t))^Ty'(t) = \nabla H(y(t))^TJ\nabla H(y(t)) = 0,$$
due to the fact that $J$ is skew-symmetric. For isolated mechanical systems, the Hamiltonian has the physical meaning of the total energy of the system and, therefore, it is of interest to derive methods which are able to preserve this property in the discrete solution. For the continuous problem, it can be seen that the {\em symplecticity} of the map implies the property of energy conservation of the given system, so that a relevant line of investigation, concerning the efficient numerical solution of such problems, has been that of devising {\em symplectic methods}, namely methods for which the discrete map inherits the property of symplecticity (see, e.g., \cite{FK85,SS88,Suris89}).
In particular, in \cite{SS88} the existence of infinitely many symplectic Runge-Kutta methods was proved, and an algebraic criterion for symplectic Runge-Kutta methods was provided.

Nevertheless, unless the continuous case, in the discrete setting the symplecticity of the map doesn't imply energy-conservation (see also \cite{BIT11}), so that a different line of investigation has been that of looking for energy-conserving methods. One of the first approaches along this line is represented by {\em discrete gradient methods} \cite{Go96,McLQR99}, which are based upon the definition of a discrete counterpart of the gradient operator, so that energy conservation for the numerical solution is guaranteed at each step and for any choice of the integration step-size. A different approach is based on the concept of {\em time finite element methods}, which has led to the definition of energy-conserving Runge-Kutta methods \cite{BS00,Bo97,TC07,TS12}, based on a local Galerkin approximations for the equation.  A partially related approach is given by {\em discrete line integral methods} \cite{IP07,IP08,IT09}, where the key idea is to exploit the relation between the method itself and the discrete line integral, i.e., the discrete counterpart of the line integral in conservative vector fields. This, in turn, allows exact conservation for polynomial Hamiltonians of arbitrarily high-degree, resulting in the class of methods later named {\em Hamiltonian Boundary Value Methods (HBVMs)}, which have been developed in a series of papers \cite{BIT09,BIT09_1,BIT10,BIS10,BIT11,BIT12,BIT12_1,BI12_1,BI12_2} (we refer to \cite{BI13} for a systematic presentation of this approach). Another approach, strictly related to the latter one, is given by the {\em averaged vector field method} \cite{QMcL08} and its generalizations \cite{Ha10}, which have been also analysed in the framework of B-series \cite{CMcLOQ10,HZ13,COS13} (i.e., methods admitting a Taylor expansion with respect to the step-size). In particular, the close connection between the limit formulae of HBVMs and the methods described in \cite{Ha10} has been thoroughly analyzed in \cite{BIT10}.

For sake of completeness, we also mention that attempts aiming to obtain methods that, in a weaker sense, have both the property of symplecticity and energy-conservation have been also considered (see, e.g., \cite{KMO99,BIT12_3,WXL13}).

Sometimes, the dynamical system defined by (\ref{ham}) has additional invariants, besides the Hamiltonian. It is therefore interesting to devise methods which are able to preserve all of them in the discrete solution. The approach based on the discrete line integrals, which HBVMs rely on, has been then used to cope with this problem, leading to the class of {\em Line Integral Methods (LIMs)} which are able to preserve any number of invariants for general conservative problems \cite{BI12} (see also \cite{BI13}).  In this paper, we consider a different generalization of HBVMs, still based on the concept of discrete line integral, which is able to provide multiple invariants conserving methods, which are more efficient than LIMs,  when the problem is in the form (\ref{ham}). 
For sake of completeness, we mention that a multiple invariants conserving version of discrete gradients is mentioned in \cite{McLQR99} (though without providing any example) and an example of such methods is given in \cite{QC99}, using an antisymmetric tensor taking discrete gradients of all integrals to be preserved as input. Additional multiple invariants conserving methods, obtained by using discrete gradients, are defined in \cite{DOY11}.

With this premise, the paper is organized as follows: in Section~\ref{1} we recall the basic facts about HBVMs; in Section~\ref{2} we define their multiple invariants conserving extension; in Section~\ref{3} we provide numerical tests for the new presented methods; finally, in Section~\ref{4} we give some conclusions.

\section{HBVMs}\label{1}
Let us consider a polynomial approximation to the solution of (\ref{ham}), over the interval $[0,h]$, in the form
\begin{equation}\label{sigma1}
\sigma'(ch) = \sum_{j=0}^{s-1} P_j(c) \gamma_j(\sigma), \qquad c\in[0,1],
\end{equation}
where $\{P_j\}_{j\ge0}$ is the family of Legendre polynomials, shifted and scaled in order to be orthonormal on the interval $[0,1]$,
\begin{equation}\label{leg}
\deg{P_j}=j, \qquad \int_0^1P_j(x)P_j(x)\dd x = \delta_{ij}, \qquad \forall i,j\ge0.
\end{equation}
By imposing the initial condition $\sigma(0)=y_0$, and setting $y_1\equiv \sigma(h)\approx y(h),$ the coefficients $\gamma_j(\sigma)$ are determined by imposing the conservation of energy at $t=h$. This implies that
\begin{eqnarray}\nonumber
0&=&H(y_1)-H(y_0) ~=~H(\sigma(h))-H(\sigma(0)) ~=~ \int_0^h \nabla H(\sigma(t))^T\sigma'(t)\dd t\\
 &=& h\int_0^1 \nabla H(\sigma(\tau h))^T\sigma'(\tau h)\dd\tau. \label{int0}\end{eqnarray}
By taking into account of (\ref{sigma1}), one then requires \cite{BIT10}:
\begin{equation}\label{zero}
\sum_{j=0}^{s-1} \left[\int_0^1 P_j(\tau) \nabla H(\sigma(\tau h))\dd\tau\right]^T\gamma_j(\sigma) = 0,
\end{equation}
which holds true, provided that
\begin{equation}\label{gammaj}
\gamma_j(\sigma) = \eta_jJ\int_0^1P_j(\tau) \nabla H(\sigma(\tau h))\dd\tau, \qquad j=0,\dots,s-1,
\end{equation}
where $\eta_0,\dots,\eta_{s-1}$ are arbitrary constants. HBVMs are then obtained by setting $$\eta_j=1,\qquad j=0,\dots,s-1,$$ resulting in an approximation of order $2s$ to $y(h)$ \cite{BIT10,BIT12}:
$$\sigma(h)-y(h) = O(h^{2s+1}).$$
In particular, by considering the orthonormality of the polynomial basis, one obtains that
$$y_1\equiv \sigma(h) = y_0 + h\gamma_0(\sigma) = y_0 + \int_0^h J\nabla H(\sigma(t))\dd t.$$
This latter expression clearly shows that this polynomial approximation generalizes that defined in \cite{QMcL08}. However, the resulting polynomial approximation, given by
\begin{equation}\label{sigma}
\sigma(ch) = y_0 + h\sum_{j=0}^{s-1} \int_0^c P_j(x)\dd x\, \gamma_j(\sigma), \qquad c\in[0,1],
\end{equation} provides an effective numerical method only when the integrals appearing in (\ref{gammaj}) are conveniently approximated by means of a quadrature formula. If this latter formula is defined at the $k$ Gauss-Legendre points in $[0,1]$,
\begin{equation}\label{ci} 0<c_1<\dots<c_k<1,\end{equation} (i.e., $P_k(c_i)=0$, $i=1,\dots,k$) and corresponding quadrature weights
\begin{equation}\label{bi} b_1, \dots, b_k >0, \end{equation}
one then obtains a HBVM$(k,s)$ method which can be cast as a $k$-stage Runge-Kutta method, with abscissae (\ref{ci}), weights (\ref{bi}), and Butcher matrix given by
\begin{equation}\label{Butmat}
A = \II_s\PP_s^T\Omega,
\end{equation}
where
$$\PP_s = \left( P_{j-1}(c_i)\right), ~\II_s = \left( \int_0^{c_i} P_{j-1}(x)\dd x\right) \,\in\RR^{k\times s}, \qquad \Omega = \diag(b_1,\dots,b_k).$$ The corresponding polynomial approximation is then given by \cite{BIT10,BIT12}
\begin{eqnarray}\nonumber
u(ch) &=& y_0 + h\sum_{j=0}^{s-1} \int_0^cP_j(x)\dd x\left[\sum_{\ell=1}^k b_\ell P_j(c_\ell)J\nabla H(u_\ell)\right]\\
 &\equiv& y_0 + h\sum_{j=0}^{s-1} \int_0^cP_j(x)\dd x\, \hgamma_j, \qquad c\in[0,1],\label{uc} \end{eqnarray} where
\begin{equation} \label{ul}
u_\ell \equiv u(c_\ell h), \qquad \ell=1,\dots,k,
\end{equation}
are nothing but the stages of the Runge-Kutta method. It can be proved that \cite{BIT09,BIT10,BIT12}, for all $k\ge s$, a HBVM$(k,s)$ method:
\begin{itemize}
\item has order $2s$;

\item is symmetric;

\item when $k=s$ it reduces to the $s$-stage Gauss-Legendre method;

\item is energy-conserving for all polynomial Hamiltonians of degree not larger than $2k/s$. Differently, the error in the Hamiltonian is $O(h^{2k+1})$, when $H$ is suitably regular.

\end{itemize}
From the last point, a {\em practical} conservation of the Hamiltonian follows, also considering that the computational complexity of the method is $s$, independently of $k$. Indeed, by reformulating the discrete problem generated by the method in terms of the $s$ unknown coefficients $\{\hgamma_j\}$ appearing in (\ref{uc}), one obtains the discrete problem \cite{BIT11}
$$\bfgamma = \PP_s^T\Omega\otimes J\,\nabla H\left( \bfe\otimes y_0 + h\II_s\otimes I\,\bfgamma\right),$$ where
$$\bfe = \pmatrix{c} 1\\ \vdots \\1\endpmatrix \in\RR^k, \qquad \bfgamma = \pmatrix{c}\hgamma_0\\ \vdots\\ \hgamma_{s-1}\endpmatrix,$$ which has (block)-size $s$, independently of $k$.

It is worth mentioning that, because of the existing relations between the integrals of the polynomials $\{P_j\}$ and the polynomials themselves, matrix (\ref{Butmat}) can be also written as
\begin{equation}\label{Wtransf}
A = \PP_{s+1} \hat{X}_s \PP_s^T\Omega,
\end{equation}
where $$ \hat{X}_s = \pmatrix{cccc}
\frac{1}2 &-\xi_1\\
\xi_1      &0   &\ddots\\
              &\ddots &\ddots &-\xi_{s-1}\\
              &           &\xi_{s-1} &0\\
                            \hline
              &           &               &\xi_s\endpmatrix \equiv \pmatrix{c} X_s\\\hline  0\,\dots\,0~\xi_s\endpmatrix,
$$ with $$\xi_i = \left( 2 \sqrt{4i^2-1}\right)^{-1}, \qquad i=1,\dots,s.$$ By considering that, for $k=s$,
$$\PP_{s+1} = \pmatrix{cc} \PP_s & 0\endpmatrix, \qquad \PP_s^T\Omega = \PP_s^{-1},$$ one then sees that (\ref{Wtransf}) can be regarded as a generalization of the $W$-transformation for collocation methods, as defined by Hairer and Wanner \cite[page\,79]{HaWa}.

\section{Multiple invariants conserving HBVMs}\label{2}
We now use again the approach based on line integrals, to define a multiple invariants conserving version of HBVM$(k,s)$ methods. Though the basic idea is similar to that used in \cite{BI12}, nevertheless, the obtained methods are definitely different from those described in that reference:  
the similarity between the two classes of methods stems from the use of the same, straightforward, methodological tool, given by discrete line integrals \cite{IP07}.  

In more details, we now use the fact that energy conservation is gained, with $\gamma_j(\sigma)$ in the form of (\ref{gammaj}), whichever  $\eta_j$ is. Assume then that
\begin{equation}\label{nuinv}
L:\RR^{2m}\rightarrow \RR^\nu
\end{equation}
is a set of $\nu$ (functionally independent) smooth invariants for the dynamical system (\ref{ham}), besides the Hamiltonian $H$. Consequently, one has
\begin{equation}\label{nuinv1}
\nabla L(y)^T J\nabla H(y) = 0\in\RR^\nu, \qquad \forall y, \end{equation}
where $\nabla L(y)^T$ is the Jacobian  matrix of $L$. We will now extend the approach described in the previous section, in order to impose their conservation. For sake of simplicity, we shall at first define a polynomial approximation $\sigma\in\Pi_s$  (where  $s>\nu$), at a continuous level (i.e., similar to (\ref{sigma1})--(\ref{sigma})), then passing to define a fully discrete approximation $u\in\Pi_s$.  Clearly, by setting $\sigma$ in the form (\ref{sigma1}), we gain energy-conservation by repeating similar steps as done until (\ref{sigma}). The difference, in such a case, is obtained by setting
\begin{eqnarray}\nonumber
\eta_j&=&1,\qquad j=0,\dots,s-\nu-1,\\ \label{etaj}\\
\eta_j &=& \left[1-h^{2(s-1-j)}\aa_j\right], \qquad j=s-\nu,\dots,s-1,\nonumber
\end{eqnarray}
with the coefficients $\{\aa_j\}$ determined in order to obtain the conservation of the $\nu$ additional invariants (\ref{nuinv})-(\ref{nuinv1}), even though, in principle, any subset of the $\nu$ coefficients $\{\eta_j\}$ could be used for this purpose.
By setting, as before, the new approximation
$$y_1\equiv \sigma(h)\approx y(h),$$ from (\ref{sigma1}), (\ref{gammaj}), and (\ref{etaj}) one then obtains, by requiring conservation of all invariants,
\begin{eqnarray*}
0 &=& L(y_1)-L(y_0) ~=~ L(\sigma(h))-L(\sigma(0)) ~=~ \int_0^h \nabla L(\sigma(t))^T\sigma'(t)\dd t \\
&=& h\int_0^1\nabla L(\sigma(\tau h))^T\sigma'(\tau h)\dd\tau ~=~ h\sum_{j=0}^{s-1} \left[\int_0^1P_j(\tau)\nabla L(\sigma(\tau h))\dd\tau\right]^T\gamma_j(\sigma)\dd\tau \\
&\equiv& h\left[\sum_{j=0}^{s-1}\phi_j(\sigma)^T\bgamma_j(\sigma) ~-~\sum_{j=s-\nu}^{s-1} h^{2(s-1-j)}\aa_j\phi_j(\sigma)^T\bgamma_j(\sigma)\right],
\end{eqnarray*}
where (see (\ref{gammaj}) and (\ref{etaj})), for all $j\geq 0$:
\begin{eqnarray}\nonumber
\phi_j(\sigma) &=& \int_0^1P_j(\tau)\nabla L(\sigma(\tau h))\dd\tau ~\in\RR^{2m\times\nu}, \\ \label{figamj}
\\[-.4cm] \nonumber
\bgamma_j(\sigma) &=& \int_0^1P_j(\tau)J\nabla H(\sigma(\tau h))\dd\tau~\in\RR^{2m}.\end{eqnarray}
Consequently, energy-conservation is ``for free'' and, moreover, the conservation of the invariants is gained provided that
\begin{equation}\label{Lcons}
\sum_{j=s-\nu}^{s-1} h^{2(s-1-j)}\aa_j\phi_j(\sigma)^T\bgamma_j(\sigma) = \sum_{j=0}^{s-1}\phi_j(\sigma)^T\bgamma_j(\sigma).
\end{equation}
By defining the matrix
\begin{equation}\label{Gam}
\Gamma(\sigma) = \left[\begin{array}{ccc} h^{2(\nu-1)}\phi_{s-\nu}(\sigma)^T\bgamma_{s-\nu}(\sigma), & \dots ~,& h^0\phi_{s-1}(\sigma)^T\bgamma_{s-1}(\sigma)\end{array}\right]\in\RR^{\nu\times\nu}\end{equation}
and the vectors
$$\balfa = \pmatrix{c} \aa_{s-\nu}\\ \vdots\\ \aa_{s-1}\endpmatrix, \quad \bfb(\sigma) = \sum_{j=0}^{s-1}\phi_j(\sigma)^T\bgamma_j(\sigma) \quad\in\RR^\nu,$$ equation (\ref{Lcons}) can be recast in vector form as
\begin{equation}\label{Lcons1}
\Gamma(\sigma) \balfa = \bfb(\sigma).
\end{equation}
The following results then hold true.

\begin{lem}\label{psij} Let $\psi:\,[0,h]\rightarrow V$, with $V$ a vector space, admit a Taylor expansion at 0. Then, for all $j\ge0$:
$$\int_0^1P_j(\tau)\psi(\tau h)\dd\tau = O(h^j).$$
\end{lem}
\proof By taking into account (\ref{leg}), one obtains:
\begin{eqnarray*}
\int_0^1P_j(\tau)\psi(\tau h)\dd\tau &=& \int_0^1 P_j(\tau) \sum_{n\ge0} \frac{\psi^{(n)}(0)}{n!} \tau^nh^n\dd\tau~=~\sum_{n\ge0} \frac{\psi^{(n)}(0)}{n!} h^n \int_0^1P_j(\tau)\tau^n\dd\tau\\
&=& \sum_{n\ge j} \frac{\psi^{(n)}(0)}{n!} h^n \int_0^1P_j(\tau)\tau^n\dd\tau ~=~O(h^j).\,\QED
\end{eqnarray*}

\begin{lem}\label{series}
If $H$ is suitably regular, then the right-hand side of problem (\ref{ham}) can be expanded as
$$J\nabla H(y(ch)) = \sum_{j\ge0} P_j(c) \bgamma_j(y), \qquad c\in[0,1],$$
where $\bgamma(y)$ is defined according to (\ref{figamj}).
\end{lem}
\proof See \cite{BIT12_1}.\,\QED

\begin{lem}\label{bfb} With reference to (\ref{Lcons1}), one has: $\bfb(\sigma) = O(h^{2s})$.
\end{lem}
\proof From (\ref{nuinv1}) and (\ref{figamj}) one obtains:
$$\sum_{j\ge0} \phi_j(\sigma)^T\bgamma_j(\sigma) = 0.$$ Consequently, by virtue of Lemma~\ref{psij},
$$\bfb(\sigma) = \sum_{j=0}^{s-1} \phi_j(\sigma)^T\bgamma_j(\sigma) = -\sum_{j\ge s}\phi_j(\sigma)^T\bgamma_j(\sigma) = O(h^{2s}).\,\QED$$

\begin{lem}\label{Gamma} Matrix $\Gamma(\sigma)$ in (\ref{Lcons1}) has  $O(h^{2s-2})$ entries.
\end{lem}
\proof The proof follows immediately from (\ref{Gam}) and Lemma~\ref{psij}.\,\QED\bigskip

In order to simplify the subsequent arguments, we make the following assumption on matrix $\Gamma(\sigma)$:\footnote{Actually, it would suffice the system (\ref{Lcons1}) to be consistent, but the arguments would become more involved.}

\begin{ass}\label{ass1} Matrix $\Gamma(\sigma)$ is nonsingular.\end{ass}

The following result then easily follows from Lemmas~\ref{bfb} and \ref{Gamma}.

\begin{theo}\label{alfa2} Under Assumption~\ref{ass1}, the vector $\balfa$ in (\ref{Lcons1}) has $O(h^2)$ entries.\end{theo}

\begin{rem}\label{nulessrem}
We observe that, in order for Theorem~\ref{alfa2} to hold, it is {\em necessary} that
\begin{equation}\label{nuless}s>\nu.\end{equation}
 In fact, when $s=\nu$,   from (\ref{Lcons}) one obtains that the products $h^{2(s-1-j)}\aa_j$ are all equal to 1 and then, from (\ref{etaj}) it follows that $\eta_j=0$,  $j=0,\dots, s-1$.  Consequently, in the sequel we shall assume that \begin{equation}\label{eta0eq1}\eta_0=1.\end{equation}
\end{rem}

We can now state the following result.

\begin{cor}\label{ord2s} Under Assumption~\ref{ass1}, the method conserves all the invariants. Moreover, $\sigma(h)-y(h)=O(h^{2s+1})$.\end{cor}

\proof The first part of the proof follows from the definition of the method. The second part of the
proof strictly follows the technique used in \cite{BIT12_1}. Let then
$y(t;\omega,z)$ be the solution of problem (\ref{ham}) satisfying the initial condition $y(\omega)=z$.
Moreover, let $\Phi(t,\tau)$ be the corresponding fundamental matrix solution of the associated variational problem. Consequently, from Lemmas~\ref{psij} and \ref{series}, Theorem~\ref{alfa2}, and from (\ref{etaj}), one obtains:
\begin{eqnarray*}
\lefteqn{\sigma(h) - y(h)~=~y(h;h,\sigma(h))-y(h;0,\sigma(0)) ~=~ \int_0^h \frac{\dd}{\dd t} y(h;t,\sigma(t))\dd t}\\
&=& \int_0^h \left[ \frac{\partial}{\partial \omega}\left.y(h;\omega,\sigma(t))\right|_{\omega=t} +
\frac{\partial}{\partial z}\left.y(h;t,z)\right|_{z=\sigma(t)}\sigma'(t)\right]\dd t\\
&=& \int_0^h \Phi(h,t)\left[ -J\nabla H(\sigma(t)) + \sigma'(t)\right] \dd t\\
&=& h\int_0^1\Phi(h,\tau h)\left[ -\sum_{j\ge0} P_j(\tau)\bgamma_j(\sigma) +\sum_{j=0}^{s-1} P_j(\tau)\gamma_j(\sigma)\right]\dd\tau\\
&=&-h\int_0^1\Phi(h,\tau h)\left[ \sum_{j\ge s} P_j(\tau)\bgamma_j(\sigma) +\sum_{j=s-\nu}^{s-1} P_j(\tau)h^{2(s-1-j)}\aa_j\bgamma_j(\sigma)\right]\dd\tau\\
&=&-h\sum_{j\ge s} \left[\underbrace{\int_0^1P_j(\tau)\Phi(h,\tau h)\dd\tau}_{=O(h^j)}\right] \overbrace{\bgamma_j(\sigma)}^{=O(h^j)}-\\
&&h\sum_{j= s-\nu}^{s-1} \left[\underbrace{\int_0^1P_j(\tau)\Phi(h,\tau h)\dd\tau}_{=O(h^j)}\right] h^{2(s-1-j)}\overbrace{\aa_j}^{=O(h^2)}\underbrace{\bgamma_j(\sigma)}_{=O(h^j)}
~=~ O(h^{2s+1}).\,\QED
\end{eqnarray*}

\subsection{Discretization and ELIM$(r,k,s)$ methods} 
As is clear, the polynomial approximation $\sigma\in\Pi_s$ defined above doesn't yet provide a numerical method: this will be obtained once the integrals in (\ref{figamj}) are approximated by means of a suitable quadrature formula. As in the case of LIM$(r,k,s)$ methods in \cite{BI12}, for this purpose we choose the abscissae 
\begin{equation}\label{hci}
0<\hc_1<\dots<\hc_r<1,
\end{equation}
placed at the $r$ Gauss-Legendre points in [0,1], and the corresponding weights 
\begin{equation}\label{hbi}
\hb_1,\dots,\hb_r>0,
\end{equation}
besides (\ref{ci})--(\ref{bi}) previously considered. In so doing, we obtain a new polynomial approximation, say $u\in\Pi_s$, defined by replacing the integrals with the given quadrature formula, having order $2r$ or $2k$, depending on the chosen abscissae. By setting $u_\ell$ formally defined as in (\ref{ul}), and (see (\ref{hci})) $$\hu_\ell\equiv u(\hc_\ell), \qquad \ell=1,\dots,r,$$ for $0\leq j\leq s-1$ one then obtains:
\pagebreak
\begin{eqnarray}\nonumber
\hphi_j &=& \sum_{\ell=1}^r \hb_\ell P_j(\hc_\ell) \nabla L(\hu_\ell) ~\equiv~ \phi_j(u) - \Psi_j(h),\\
\label{figamjd}\\[-.4cm]
\hgamma_j &=&  \sum_{\ell=1}^k b_\ell P_j(c_\ell) J\nabla H(u_\ell) ~\equiv~ \bgamma_j(u) - \Delta_j(h),
\nonumber
\end{eqnarray}
in place of (\ref{figamj}) where, by denoting
\begin{equation}\label{mu}
\mu(j) = \lfloor \frac{2j}s\rfloor, \qquad j\in\{r,\,k\},
\end{equation}
and assuming $L$ and $H$ suitably regular,
\begin{equation}\label{Psih}
\Psi_j(h) = \left\{ \begin{array}{ccc} 0, &\quad &\mbox{if}\quad L\in\Pi_{\mu(r)},\\[.4cm]
O(h^{2r-j}), &&\mbox{otherwise,}\end{array}\right.
\end{equation}
and
\begin{equation}\label{Deltah}
\Delta_j(h) = \left\{ \begin{array}{ccc} 0, &\quad &\mbox{if}\quad H\in\Pi_{\mu(k)},\\[.4cm]
O(h^{2k-j}), &&\mbox{otherwise.}\end{array}\right.
\end{equation}
\begin{rem} 
Actually, for any invariant in (\ref{figamjd}) one could use a different quadrature formula, depending on the required accuracy. Nevertheless, for sake of brevity, we shall hereafter consider only the use of two (possibly) different quadratures: (\ref{ci})--(\ref{bi}) for the $\{\hgamma_j\}$, and (\ref{hci})-(\ref{hbi}) for the $\{\hphi_j\}$. However, the generalization is straightforward. 
\end{rem}
Setting by $\heta_j$ and $\haa_j$, respectively, the discrete approximations to (\ref{etaj}), the new polynomial approximation is then given by
\begin{eqnarray}\label{ul1}
u(ch) &=& y_0 + h\sum_{j=0}^{s-1} \int_0^cP_j(x)\dd x \, \heta_j \hgamma_j\\
 &\equiv& y_0 + h\left[\sum_{j=0}^{s-1} \int_0^cP_j(x)\dd x\, \hgamma_j -  \sum_{j=s-\nu}^{s-1}  \int_0^cP_j(x)\dd x\, h^{2(s-1-j)}\hat\aa_j \hgamma_j\right], \nonumber
\end{eqnarray}
with the scalars $\haa_j$ satisfying the equation (compare with (\ref{Lcons})):
\begin{equation}\label{Lconsd}
\sum_{j=s-\nu}^{s-1} h^{2(s-1-j)}\haa_j\hphi_j^T\hgamma_j = \sum_{j=0}^{s-1}\hphi_j^T\hgamma_j.
\end{equation}
Similarly as previously done in (\ref{Gam})--(\ref{Lcons1}), by defining the matrix
\begin{equation}\label{Gamd}
\hat\Gamma = \left[\begin{array}{ccc} h^{2(\nu-1)}\hphi_{s-\nu}^T\hgamma_{s-\nu}, & \dots~, & h^0\hphi_{s-1}^T\hgamma_{s-1}\end{array}\right]\in\RR^{\nu\times\nu}\end{equation}
and the vectors
$$\hat{\balfa} = \pmatrix{c} \haa_{s-\nu}\\ \vdots\\ \haa_{s-1}\endpmatrix, \quad \hat{\bfb} = \sum_{j=0}^{s-1}\hphi_j^T\hgamma_j \quad\in\RR^\nu,$$ equation (\ref{Lconsd}) can be recast in vector form as
\begin{equation}\label{Lcons1d}
\hat{\Gamma} \hat{\balfa} = \hat{\bfb}.
\end{equation}
Since the number of the additional invariants (\ref{nuinv1}) has to satisfy (\ref{nuless}) (and, then, (\ref{eta0eq1}) holds true), similarly as in the case of HBVM$(k,s)$, the new approximation is  given by
\begin{equation}\label{y1uh}
y_1\equiv u(h) = y_0 + h\hgamma_0 = y_0 + h\sum_{\ell=1}^k b_\ell J\nabla H(u_\ell).
\end{equation}
\begin{defi}\label{elimdef} 
 We shall denote by ELIM$(r,k,s)$ (Enhanced LIM$(r,k,s)$) the methods defined by (\ref{ul1})--(\ref{y1uh}). In particular, for similarity with the GHBVM$(k,s)\,\equiv\,$LIM$(k,k,s)$ methods in \cite{BI12}, when $r=k$ we shall speak about an EHBVM$(k,s)$ (Enhanced HBVM$(k,s)$) method. 
 \end{defi}

The following results then easily follow, providing a discrete counterpart of Theorem~\ref{alfa2}.

\begin{theo}\label{alfa2d} 
Under Assumption~\ref{ass1}, for all ~$r,k\ge s$~ matrix $\hat{\Gamma}$ is nonsingular, for all  sufficiently small step-sizes $h$, and the vector $\hat{\balfa}$ has $O(h^2)$  entries.\end{theo}

We can now state the following results, concerning the order of accuracy of the discrete solution, as well as of the invariants, provided by ELIM$(r,k,s)$ methods.

\begin{theo}\label{exactH} 
Assume that the Hamiltonian function defining problem (\ref{ham}) is a polynomial of degree less than or equal to $\mu(k)$ as defined in (\ref{mu}). Then, ELIM$(r,k,s)$ method is energy-conserving, provided that $r\ge s$.  Differently, for all general  and suitably regular $H$, one obtains $$H(y_1)-H(y_0) = O(h^{2k+1}), \qquad \forall k\ge s,$$ provided that $r\ge s$.
\end{theo}
\proof
One has:
\begin{eqnarray*}
\lefteqn{H(y_1)-H(y_0) ~=~ H(u(h))-H(u(0)) ~=~\int_0^h \nabla H(u(t))^Tu'(t)\dd t}\\
&=& h\int_0^1  \nabla H(u(\tau h))^Tu'(\tau h)\dd \tau ~=~ h \int_0^1 \nabla H(u(\tau h))^T \sum_{j=0}^{s-1}P_j(\tau) \heta_j\hgamma_j\dd\tau\\
&=& h \sum_{j=0}^{s-1} \left[\int_0^1 \nabla H(u(\tau h))P_j(\tau)\dd\tau\right]^T \heta_j\hgamma_j
~=~ h\sum_{j=0}^{s-1} \heta_j \bgamma_j(u)^TJ\hgamma_j ~=~(*).
\end{eqnarray*}
The first part of the proof easily follows from the fact that, if $H\in\Pi_{\mu(k)}$, then
$$\hgamma_j = \bgamma_j(u), \qquad j=0,\dots,s-1,$$ so that $(*)=0$, since $J$ is skew-symmetric.
In general, assuming that $H$ is suitably regular,  one has (see (\ref{figamjd}) and (\ref{Deltah})):
$$
(*) = h\sum_{j=0}^{s-1} \heta_j \bgamma_j(u)^TJ(\bgamma_j(u) - \Delta_j(h))
 ~=~ -h\sum_{j=0}^{s-1} \overbrace{\heta_j}^{=O(1)} \underbrace{\bgamma_j(u)^T}_{=O(h^j)}J\overbrace{\Delta_j(h)}^{=O(h^{2k-j})} ~=~ O(h^{2k+1}).\,\QED
$$

Using similar arguments, by means of (\ref{Psih}) it is possible to prove the following result.

\begin{theo}\label{exactL} 
Assume that the invariants (\ref{nuinv1}) of problem (\ref{ham}) are polynomials of degree less than or equal to $\mu(r)$ as defined in (\ref{mu}). Then, EHBVM$(r,k,s)$ method is invariants-conserving, provided that $k\ge s$. For all general and suitably regular $L$, one obtains $$L(y_1)-L(y_0) = O(h^{2r+1}), \qquad \forall r\ge s,$$provided that $k\ge s$.
\end{theo}

Next result concerns the order of accuracy of the numerical solution.

\begin{theo}\label{ordu1}  
Assuming that both $H$ and $L$ are suitably regular, for all ~$r,k\ge s$~ the numerical solution generated by a ELIM$(r,k,s)$ method satisfies $$y_1-y(h)=O(h^{2s+1}).$$ That is, the method has order $2s$.\end{theo}

\proof The proof proceeds in a similar way as that of Corollary~\ref{ord2s}. By using the same notation in that corollary, one has:
\begin{eqnarray*}
\lefteqn{y_1-y(h) ~=~ y(h;h,u(h))-y(h,0,u(0)) ~=~ \int_0^h \frac{\dd}{\dd t} y(h;t,u(t))\dd t}\\
&=& \int_0^h \left[ \frac{\partial}{\partial \omega}\left.y(h;\omega,u(t))\right|_{\omega=t} +
\frac{\partial}{\partial z}\left.y(h;t,z)\right|_{z=u(t)}u'(t)\right]\dd t\\
&=& \int_0^h \Phi(h,t)\left[ -J\nabla H(u(t)) + u'(t)\right] \dd t\\
&=& h\int_0^1\Phi(h,\tau h)\left[ -\sum_{j\ge0} P_j(\tau)\bgamma_j(u) +\sum_{j=0}^{s-1} P_j(\tau)\heta_j\hgamma_j\right]\dd\tau\\
&=& h\int_0^1\Phi(h,\tau h)\left[ -\sum_{j\ge0} P_j(\tau)\bgamma_j(u) +\sum_{j=0}^{s-1} P_j(\tau)\heta_j(\bgamma_j(u)-\Delta_j(h))\right]\dd\tau\\
&=&h\int_0^1\Phi(h,\tau h)\left[-\sum_{j=0}^{s-1} P_j(\tau)\Delta_j(h) - \sum_{j\ge s} P_j(\tau)\bgamma_j(u) -\sum_{j=s-\nu}^{s-1} P_j(\tau)h^{2(s-1-j)}\haa_j\bgamma_j(u)\right]\dd\tau\\
&=&h\sum_{j=0}^{s-1} \left[\underbrace{-\int_0^1P_j(\tau)\Phi(h,\tau h)\dd\tau}_{=O(h^j)}\right] \overbrace{\Delta_j(h)}^{=O(h^{2k-j})}-
h\sum_{j\ge s} \left[\underbrace{\int_0^1P_j(\tau)\Phi(h,\tau h)\dd\tau}_{=O(h^j)}\right] \overbrace{\bgamma_j(u)}^{=O(h^j)}-\\
&&h\sum_{j= s-\nu}^{s-1} \left[\underbrace{\int_0^1P_j(\tau)\Phi(h,\tau h)\dd\tau}_{=O(h^j)}\right] h^{2(s-1-j)}\overbrace{\haa_j}^{=O(h^2)}\underbrace{\bgamma_j(u)}_{=O(h^j)}\\
&=& O(h^{2k+1}) + O(h^{2s+1}) + O(h^{2s+1}) ~=~ O(h^{2s+1}).\,\QED
\end{eqnarray*}

For sake of completeness, we also mention the following result, whose proof is straightforward (see, e.g., \cite{BIT09,BI13}).

\begin{theo}  
Provided that the abscissae (\ref{ci}) and (\ref{hci}) are symmetrically distributed in the interval [0,1], the ELIM$(r,k,s)$ method is symmetric. 
\end{theo}

\subsection{Runge-Kutta type formulation of ELIM$(r,k,s)$ methods}
Though the method (\ref{ul1})--(\ref{y1uh}) is not strictly a Runge-Kutta method, nevertheless, it admits a  Runge-Kutta type formulation which is quite useful to represent it. In more details, we already saw that a HBVM$(k,s)$ methods is a $k$-stage Runge-Kutta method defined by the following Butcher tableau (see (\ref{Butmat}))
$$\begin{array}{c|c}
\bfc & \II_s\PP_s^T\Omega\\
\hline &\bfb^T\end{array},$$
where, as usual, $\bfc$ is the vector of the abscissae and $\bfb$ is the vector of the weights. Moreover, we recall that the only formal difference between a HBVM$(k,s)$ method and   
an ELIM$(r,k,s)$ method consists in the coefficients $\heta_1,\dots,\heta_{s-1}$  
which may assume values different from 1 (indeed, $\heta_0=1$, as stated in (\ref{eta0eq1})). Consequently, by introducing the diagonal matrix
$$\Sigma_s = \diag(1,\heta_1,\dots,\heta_{s-1}),$$ one obtains the following Runge-Kutta type formulation of an ELIM$(r,k,s)$ method:
$$\begin{array}{c|c}
\bfc & \II_s\Sigma_s\PP_s^T\Omega\\
\hline &\bfb^T\end{array}.$$
As an example, HBVM(2,2) is the usual 2-stage Gauss method, whereas ELIM($r$,2,2)  
is given by {\small
$$
\begin{array}{c|cc}
\frac{1}2-\frac{\sqrt{3}}6 & \frac{1}4+(\heta_1-1)\frac{\sqrt{3}}{12} & \frac{1}4-(\heta_1+1)\frac{\sqrt{3}}{12} \\[.2cm]
\frac{1}2+\frac{\sqrt{3}}6 & \frac{1}4+(\heta_1+1)\frac{\sqrt{3}}{12} &
\frac{1}4-(\heta_1-1)\frac{\sqrt{3}}{12} \\[.2cm]
\hline
& \frac{1}2 &\frac{1}2
\end{array}
$$}
As expected, when $\heta_1=1$ one retrieves the usual 2-stage Gauss method.

\section{Numerical tests}\label{3}
We here report a few numerical tests, aimed to assess the theoretical findings, as well as to compare the Enhanced Line Integral Methods (ELIMs), here introduced, with the Line Integral Methods (LIMs) defined in \cite{BI12}. This will be done on a Hamiltonian problem possessing multiple invariants. The generated discrete problems are solved by means of a fixed-point iteration, even though the efficient implementation of both classes of methods deserves a further investigation. 

In order to compare the methods, it will be useful to consider that, for a given problem possessing $\nu$ invariants besides the Hamiltonian, one has:
\begin{equation}\label{costratio}
\frac{\mbox{cost of 1 LIM$(r_1,k_1,s)$ fixed-point iteration}}{\mbox{cost of 1 ELIM$(r_2,k_2,s)$ fixed-point iteration}} \approx \frac{k_1+(\nu+1) r_1}{k_2+\nu r_2}.
\end{equation}
Consequently, when $r_1=k_1=r_2=k_2\equiv k$, one obtains that\footnote{We recall that \cite{BI12} LIM$(k,k,s)\equiv$\,GHBVM$(k,s)$, and 
(see Definition~\ref{elimdef}) ELIM$(k,k,s)\equiv$\,EHBVM$(k,s)$.}
\begin{equation}\label{costratio1}
\frac{\mbox{cost of 1 GHBVM$(k,s)$ fixed-point iteration}}{\mbox{cost of 1 EHBVM$(k,s)$ fixed-point iteration}} \approx \frac{\nu+2}{\nu +1}.
\end{equation}

That said, the problem that we consider is the well known Kepler problem \cite{HLW06,BI12}, defined by the Hamiltonian
\begin{equation}\label{keplH}
H(q,p) = \frac{1}2\|p\|_2^2 +\frac{1}{\|q\|_2},\qquad q,p\in\RR^2.
\end{equation} 
When the initial condition is chosen as
$$(q_0^T,p_0^T) = \pmatrix{cccc} 1-\eps,&0, &0, &\sqrt{\frac{1+\eps}{1-\eps}}\endpmatrix,\qquad \eps\in[0,1),$$ 
its solution is periodic, with period $2\pi$, and is given by an ellipse of eccentricity $\eps$ in the $q$-plane.
This problem admits two additional (independent) invariants of motion, besides the Hamiltonian (\ref{keplH}), given by the
{\em angular momentum}
\begin{equation}\label{keplL1}
L_1(q,p) = q^TJ_2p, \qquad J_2 = \pmatrix{cc} 0 &1\\-1&0\endpmatrix,
\end{equation}
and the {\em Laplace-Runge-Lenz (LRL) vector}, resulting in the following conserved quantity:
\begin{equation}\label{keplL2}
L_2(q,p) = (e_1^Tp)\,L_1(q,p)-\frac{e_2^Tq}{\|q\|_2},
\end{equation}
where, as usual, $e_1,e_2\in\RR^2$ are the two unit vectors.

\begin{figure}[t]
\centerline{\includegraphics[width=12cm,height=8.5cm]{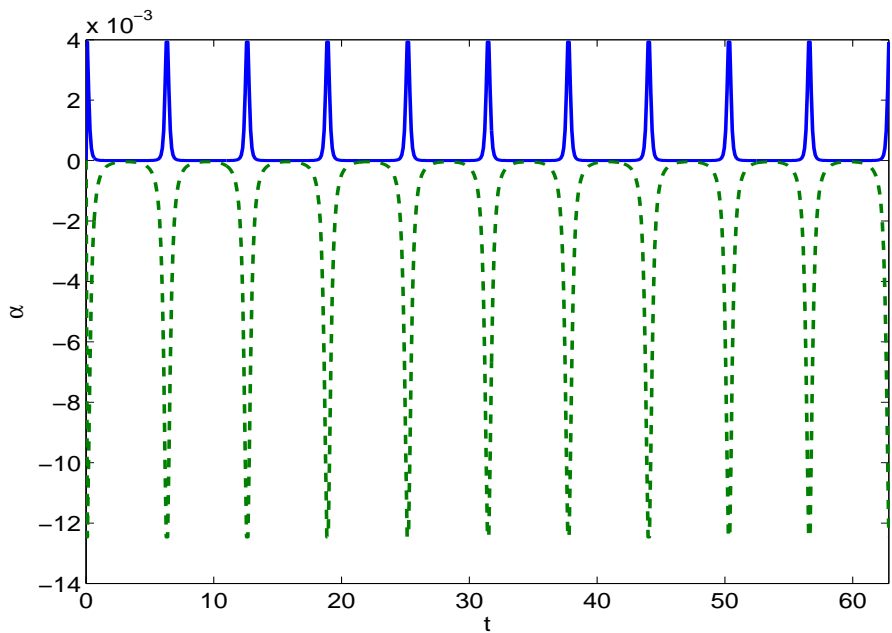}}
\caption{Components of the vector $\hat{\balfa}$ for the fully conservative EHBVM(12,3) method, $h=\pi/30$.}\label{alfa}
\end{figure}

We solve this problem, considering an eccentricity $\eps=0.6$, by using the following methods:
\begin{itemize}

\item the symplectic 3-stage Gauss method (GAUSS3);

\item the (practically) energy-conserving HBVM(12,3) method;

\item the EHBVM(12,3) method (i.e., ELIM(12,12,3)) and the GHBVM(12,3) method (i.e., LIM(12,12,3)) in \cite{BI12}, where it is imposed {\em only} the (practical) conservation of the angular momentum (\ref{keplL1}) besides the Hamiltonian (\ref{keplH});

\item the EHBVM(12,3) and GHBVM(12,3) methods as above, where it is imposed {\em both} the (practical) conservation of the angular momentum (\ref{keplL1}) and of the LRL vector (\ref{keplL2}) besides the Hamiltonian (\ref{keplH}).
\end{itemize}
In Table~\ref{erry} we list the measured errors after 10 periods, thus confirming that, according to Theorem~\ref{ordu1}, all methods are sixth-order. Moreover, in Table~\ref{alphah} we list the maximum norm for the vector $\hat{\balfa}$ defined in (\ref{Lcons1d}), over the same interval, for the EHBVM(12,3) method:
\begin{itemize}
\item by imposing {\em only} the conservation of the angular momentum besides the Hamiltonian. Here, $\alpha_h^{(1)}=\max_{n=1,\dots,\frac{T}h}\|\hat{\balfa}_n\|_\infty$;

\item by imposing {\em both} the conservation of the angular momentum and of the LRL vector besides the Hamiltonian. As before, $\alpha_h^{(2)}=\max_{n=1,\dots,\frac{T}h}\|\hat{\balfa}_n\|_\infty$.
\end{itemize}
The obtained results confirm that the entries of the vector $\hat{\balfa}$ are actually $O(h^2)$, as predicted by Theorem~\ref{alfa2d}. For sake of completeness, in Figure~\ref{alfa} we plot the two components of the vector $\hat{\balfa}$ in the second case, when a step-size $h=\pi/30$ is used: their periodic behavior, in accordance with that of the solution, is clearly evident.

In order to compare the computational costs of EHBVM(12,3) and GHBVM(12,3), in Table~\ref{iter} we also list the total number of fixed-point iterations needed for solving the discrete problems generated when computing the results listed in Table~\ref{erry}. From Table~\ref{iter}, one sees that GHBVM(12,3) requires approximately the same number of iterations as those needed by GAUSS3 and HBVM(12,3) methods (this fact was already known from \cite{BI12}), whereas EHBVM(12,3) requires some extra iteration, which increase with the number of conserved invariants. However, according to (\ref{costratio1}) one fixed-point iteration for GHBVM(12,3), when preserving $\nu$ invariants besides the Hamiltonian, costs $$\frac{2+\nu}{1+\nu}, \qquad \nu=1,2,$$ times than that of the corresponding EHBVM(12,3) method. This, in turn, shows that, for the considered problem, EHBVMs are more efficient than GHBVMs.

\begin{figure}[t]
\centerline{\includegraphics[width=12cm,height=8.5cm]{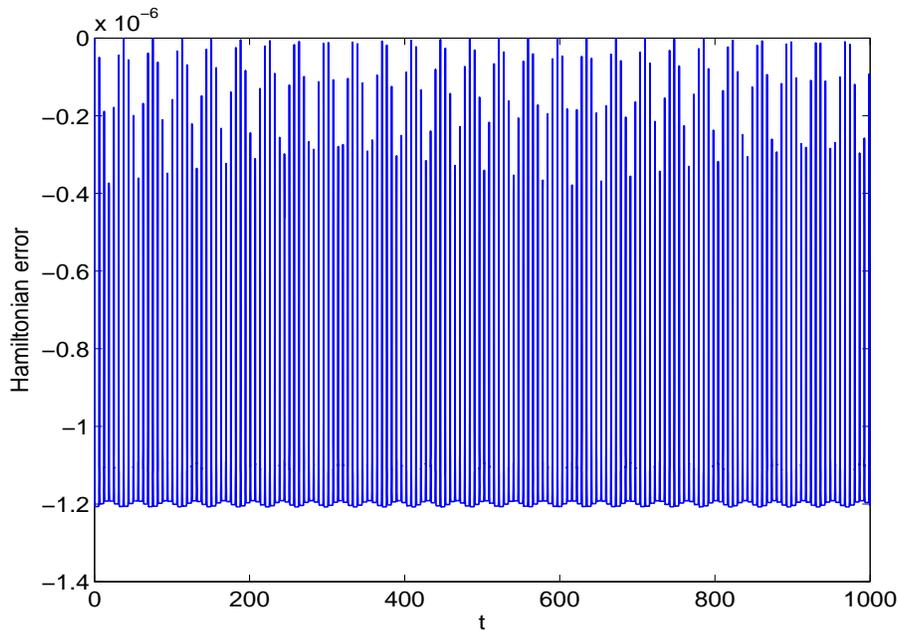}}
\caption{Hamiltonian error by using the 3-stage Gauss method, $h=0.1$.}\label{gauss_h}
\end{figure}

\begin{figure}[t]
\centerline{\includegraphics[width=12cm,height=8.5cm]{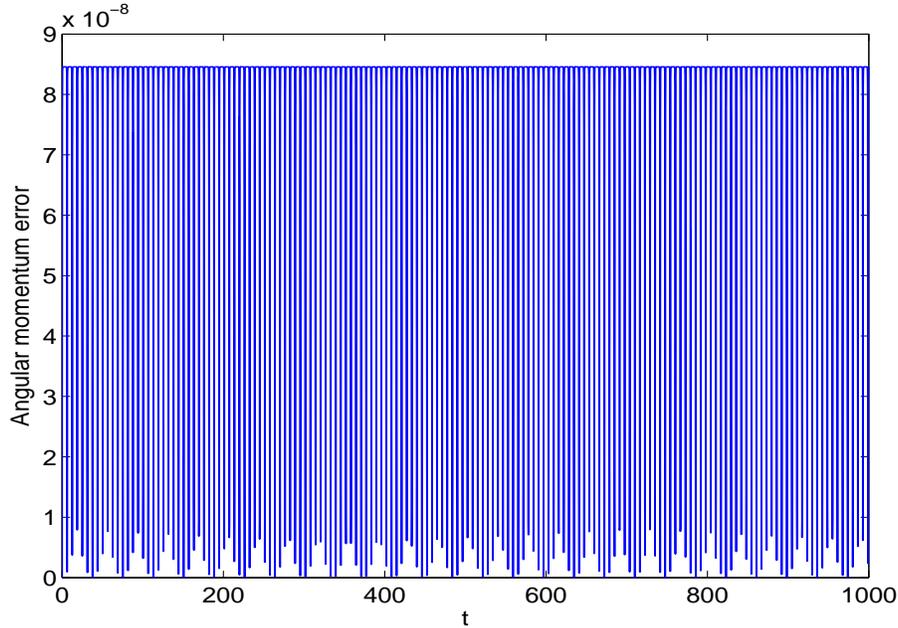}}
\caption{Angular momentum error by using the HBVM(12,3) method, $h=0.1$.}\label{hbvm_am}
\end{figure}

\begin{figure}[t]
\centerline{\includegraphics[width=12cm,height=8.5cm]{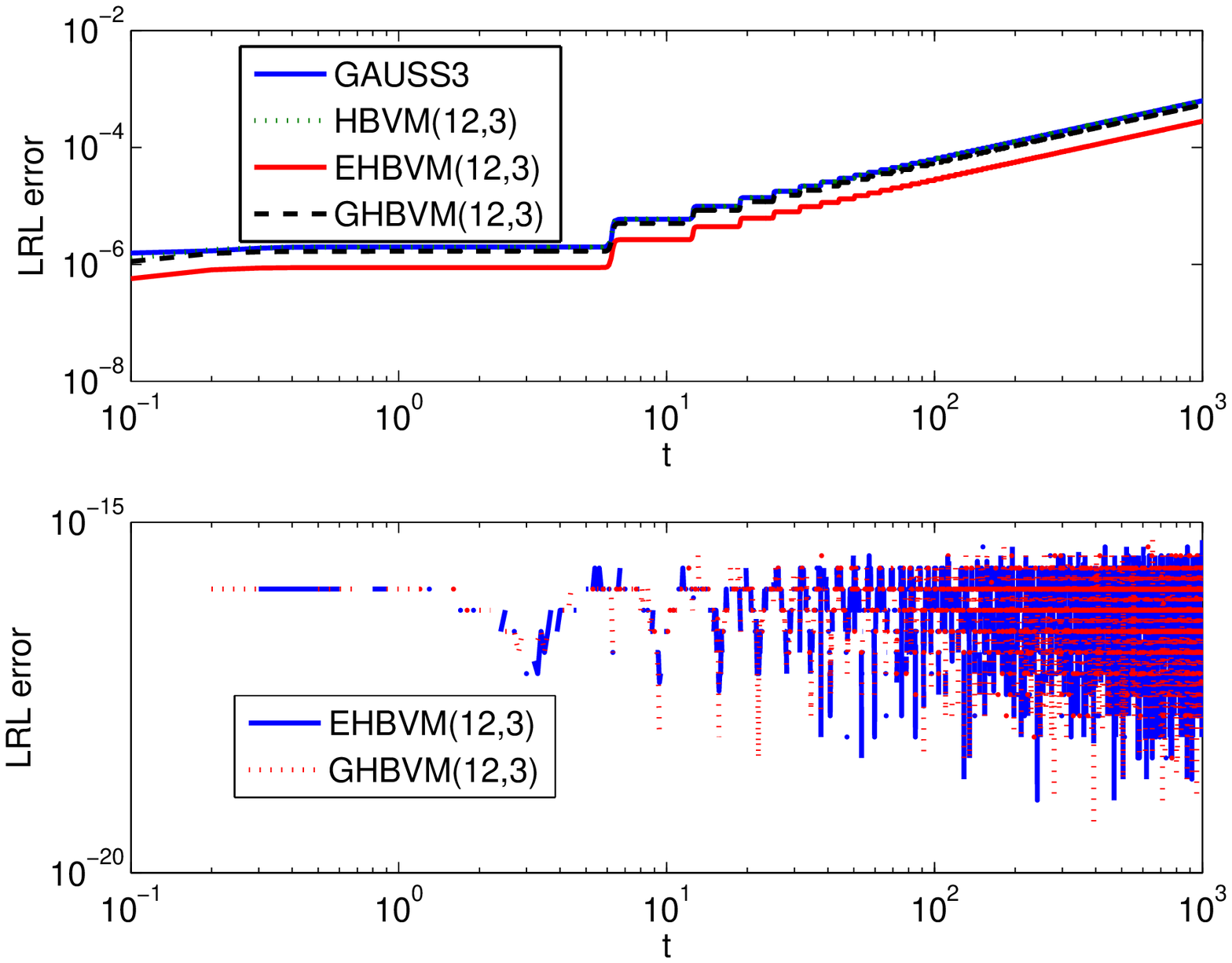}}
\caption{Error in  the LRL vector, $h=0.1$.}\label{lrl}
\end{figure}

At last, concerning the conservation of the invariants, by using a constant step-size $h=0.1$, we have solved the problem over the interval $[0,10^3]$, obtaining the following results:
\begin{itemize}

\item concerning the conservation of the Hamiltonian (\ref{keplH}), all methods are (practically) energy-conserving, except the symplectic 3-stage Gauss meth\-od. However, the Hamiltonian error turns out to be bounded, as expected, as confirmed by the plot in Figure~\ref{gauss_h};

\item concerning the conservation of the angular momentum (\ref{keplL1}), all methods conserve this invariant, except HBVM(12,3). However, the error appears to be bounded, as is shown  in Figure~\ref{hbvm_am};

\item
concerning the conservation of the LRL vector (\ref{keplL2}), all methods exhibit a drift, except EHBVM(12,3) and GHBVM(12,3), when this invariant is required to be conserved, as is shown in Figure~\ref{lrl}. In particular: the drifts of the GAUSS3 and HBVM(12,3) methods are practically the same. Both of them are slightly larger than that shown by the GHBVM(12,3) method which is, in turn, larger than that of EHBVM(12,3) method, when only the invariants (\ref{keplH}) and (\ref{keplL1}) are imposed to be conserved.\footnote{These results agree with the analysis in \cite{CLMR11}.}  

\end{itemize}

\begin{table}[t]
\caption{Error in the numerical solution for the 3-stage Gauss method ($E_G$), HBVM(12,3) ($E_H$), and EHBVM(12,3) and GHBVM(12,3) with only angular momentum conserved ($E_E^{(1)}$ and $E_G^{(1)}$, respectively), and with both angular momentum and LRL vector conserved ($E_E^{(2)}$ and $E_G^{(2)}$, respectively).} \label{erry}
\centerline{\scriptsize
\begin{tabular}{|c|cr|cr|cr|cr|cr|cr|}
\hline
$h$ & $E_G$ & order & $E_H$ & order & $E_E^{(1)}$ & order & $E_G^{(1)}$ & order & $E_E^{(2)}$ & order & $E_G^{(2)}$ & order\\
\hline
   $\pi/30~$  &   1.942e-03  &   --   & 4.587e-05 &  --  & 1.017e-05 &      --  & 4.049e-05 &  -- & 1.928e-05 & -- & 4.367e-05 &-- \\
   $\pi/60~$  &    2.817e-05 &  6.1 & 7.375e-07 & 6.0 & 1.644e-07 &  6.0 & 6.505e-07 & 6.0  & 3.052e-07 & 6.0  & 6.868e-07 & 6.0 \\
   $\pi/120$  &    4.346e-07 &  6.0 & 1.161e-08 & 6.0 & 2.591e-09 &  6.0 & 1.023e-08 &  6.0 & 4.785e-09 & 6.0  & 1.075e-08 & 6.0 \\
   $\pi/240$  &    6.771e-09 &  6.0 & 1.816e-10 & 6.0 & 4.030e-11 &  6.0 & 1.599e-10 &  6.0 & 7.509e-11 &  6.0  & 1.677e-10 & 6.0 \\
   $\pi/480$  &    1.052e-10 &  6.0 & 1.815e-12 & 6.6 & 4.718e-13 &  6.4  & 2.201e-12 & 6.2 & 1.413e-12 &  5.7  & 2.346e-12 & 6.2 \\
\hline
\end{tabular}}
\end{table}

\begin{table}[t]
\caption{Total number of fixed-point iterations for solving the discrete problems when using the 3-stage Gauss method (GAUSS3), HBVM(12,3) method, and EHBVM(12,3) and GHBVM(12,3) with only angular momentum conserved (EHBVM$_1$(12,3) and GHBVM$_1$(12,3), respectively), and with both angular momentum and LRL vector conserved (EHBVM$_2$(12,3) and GHBVM$_2$(12,3), respectively).} \label{iter}
\centerline{\scriptsize
\begin{tabular}{|c|r|r|r|r|r|r|}
\hline
$h$ & GAUSS3 & HBVM(12,3) & EHBVM$_1$(12,3) & GHBVM$_1$(12,3) & EHBVM$_2$(12,3) & GHBVM$_2$(12,3)\\
\hline
   $\pi/30~$  &    6705 &    6775 &   7256 &   6779 & 7474 & 6781 \\
   $\pi/60~$  &   11147 &  11244 & 12691 &  11247 &  13407 &11249 \\
   $\pi/120$  &   19085 & 19343 & 21664 &  19339 &  23012 &  19348\\
   $\pi/240$  &   33876 & 34752 & 37511 &  34743 & 39117 & 34753 \\
   $\pi/480$  &   61501 & 61959 &  65125 & 61967 & 68217 & 61970\\
\hline
\end{tabular}}
\end{table}

\begin{table}[t]
\caption{Quadratic convergence of the maximum norm of the vector $\hat{\balfa}$, for the Kepler problem, by using the EHBVM(12,3) method, when imposing only the angular momentum conservation ($\aa_h^{(1)}$),  and both angular momentum and LRL vector conservation ($\aa_h^{(2)}$).}\label{alphah}
\centerline{\small\begin{tabular}{|c|cc|cc|}
\hline
$h$ & $\alpha_h^{(1)}$ & order & $\alpha_h^{(2)}$ & order\\
\hline
   $\pi/30~$  &  4.530e-3 & --   & 1.246e-2 &   -- \\
   $\pi/60~$  &  1.155e-3 & 2.0 & 3.195e-3 &  2.0\\
   $\pi/120$ &  2.902e-4 & 2.0 & 8.040e-4 &   2.0\\
   $\pi/240$ &  7.265e-5 & 2.0 & 2.013e-4 &   2.0\\
   $\pi/480$ &  1.837e-5 & 2.0 & 5.055e-5 &   2.0\\
\hline
\end{tabular}}
\end{table}


\section{Conclusions}\label{4}

In this paper, we have used the technique of {\em discrete line integrals} introduced by Iavernaro and Pace \cite{IP07} to define an extension of the energy-conserving methods named HBVMs, in order to cope with the conservation of multiple invariants for Hamiltonian problems. This has resulted in an ``enhanced'' version of the {\em Line Integral Methods (LIMs)} introduced in \cite{BI12}. Consequently, we have named the new methods {\em Enhanced Line Integral Methods (ELIMs)} . The analysis of such methods has been carried out, proving that the original order of HBVMs is retained by the new methods. At last, a few numerical tests clearly confirm the theoretical findings.

\subsubsection*{Acknowledgements} This paper emerged from the visit of the first author at the Academy of Mathematics and Systems Science, Chinese Academy of Sciences, Beijing, China, in December 2012--January 2013. This has been possible because of the support provided by the Academy. The second author was supported by  the Foundation for Innovative Research Groups of the NNSFC (11021101).

The authors wish to tank an anonymous referee, for his comments and suggestions, which helped to improve the original manuscript.

\end{document}